\documentclass[12pt]{article}
\usepackage{amssymb}
\newenvironment{proof}[1][\it Proof]{\noindent \textrm{#1.} }{\hfill
\rule{0.5em}{0.5em}}
\newcommand{\bea}{\begin{eqnarray}}
\newcommand{\ena}{\end{eqnarray}}
\newcommand{\beas}{\begin{eqnarray*}}
\newcommand{\enas}{\end{eqnarray*}}
\newcommand{\beq}{\begin{equation}}
\newcommand{\enq}{\end{equation}}
\newcommand{\n}{n}
\def\qed{\hfill \mbox{\rule{0.5em}{0.5em}}}

\newcommand{\ignore}[1]{}

\newcommand{\openR}{\hbox{I\kern-.2em R}}

\title{{\bf\Large A Statistical Characterization of Regular Simplices}
}
\author{Ian Abramson and Larry Goldstein}
 \date{\empty}
 
\begin{document}
\maketitle

  \section{INTRODUCTION.}
Picture three points at the vertices of an equilateral triangle in
two dimensions, or four points at the vertices of a regular
tetrahedron in three dimensions. Thought of as scatterings of data
they wouldn't seem to reveal strong linear associations between
the coordinates. There are no clear axes of elongation in the
scatterplots, which would suggest that change in some variable is
predictable as a function of the others. In general, such
associations are usually indicated by the covariance matrix ${\bf
S}_{\bf u}$ of the set of points ${\bf u}=\{{\bf x}_1,\ldots, {\bf x}_n\}$
in $\openR^p$, which is given by
\beas
%\label{def-S}
{\bf S}_{\bf u}= \frac{1}{|{\bf u}|}\sum_{{\bf x} \in {\bf u}}
({\bf x}-\overline{\bf x}_{\bf u})({\bf x}-\overline{\bf x}_{\bf
u})',
\enas
where
 \beas
  \overline{\bf x}_{\bf u} =
\frac{1}{|{\bf u}|} \sum_{{\bf x} \in {\bf u}} {\bf x}.
\enas

The off-diagonal entries of ${\bf S}_{\bf u}$, the pairwise
covariances, tell us something about dependencies. If the
coordinate variables are independent these entries are zero.
Though the converse is false, a diagonal covariance matrix roughly
says that the coordinates are not mutually linearly predictable
from each other. Indeed, for our equilateral triangle in
$\openR^2$, tetrahedron in $\openR^3$, and the generalized
configurations in higher dimensions having equal interpoint
distances,  the covariance matrix turns out to be diagonal. In
fact, it's a scalar multiple of the identity. Furthermore, the
converse is also true: any configuration of $n=p+1$ points in $p$
dimensions whose covariance matrix is a positive multiple of the
identity are equidistant from each other. We formalize this result
in the following theorem:\\[1ex]

\noindent {\bf Theorem.} {\em Let ${\bf u}=\{{\bf x}_1,\ldots,{\bf
x}_n\}$ be a set of $n$ points in $\openR^p$, with $n = p+1 \ge
2$, and let $\sigma^2$ be an arbitrary positive number. Then the
interpoint distances of ${\bf u}$ satisfy $||{\bf x}_i-{\bf
x}_j||^2=2\sigma^2\delta_{ij}$ if and only if $n{\bf S}_{\bf
u}=\sigma^2 {\bf I}_p$.}\\[1ex]

In other words, $p+1$ points in $p$ dimensions lie at the vertices
of a regular simplex if and only if their covariance matrix is a
multiple of the identity. A proof of this statistical
characterization of regular simplices is given in section
\ref{proof}, after some preliminaries.

\section{STATISTICAL CHARACTERIZATION OF REGULAR SIMPLICES.}
\label{proof}
 The reader is assumed to be familiar with the basic
elements of linear algebra  in $\openR^p$ (linear subspaces, span,
linear dependence and independence, basis and dimension), as
treated, for example, in the text of Seber \cite{1}. For a finite
subset ${\bf u} $ of $ \openR^p$ let ${\cal V}_{\bf u} =
\mbox{span}\{{\bf x}-{\bar {\bf x}_{\bf u}}: {\bf x} \in {\bf u}
\}$.\\[1ex]

\noindent {\bf Lemma.} {\em With $\n>1$ let ${\bf u}$ be any
collection of $\n$ points in $\openR^p$ with common squared
interpoint distance $2\sigma^2>0$. Then $\mbox{dim}\left( {\cal
V}_{\bf u} \right) = \n-1$, and with
$r_{\sigma,\n}^2=\sigma^2(\n-1)/\n$ and
$s_{\sigma,\n}^2=\sigma^2/(\n(\n-1))$, the following are true for
each ${\bf x}$ in $\bf u$:
$$
||{\bf x} - {\bar {\bf x}_{\bf
u}}|| = r_{\sigma,\n}, \quad ||{\bar {\bf x}_{\bf u}} - {\bar {\bf
x}_{{\bf u} - \{\bf x\}}}|| = s_{\sigma,\n},  \quad {\bf x} - {\bar
{\bf x}_{{\bf u} - \{\bf x\}}} \perp
{\cal V}_{{\bf u} - \{{\bf x}\}}.
$$}

\begin{proof}
 We argue by induction.
The three claims are easily verified if $\n=2$. When $\n > 2$, for
every ${\bf x}$ in ${\bf u}$ the points ${{\bf u} - \{{\bf x}\}}$
are equidistant from ${\bf x}$, and by the induction hypotheses
also equidistant from their average
 ${\bar {\bf x}_{{\bf u}- {\{\bf
x}\}}}$, albeit at a smaller distance. Hence, the points of ${{\bf
u}- \{{\bf x}\}}$ lie on the intersection of two spheres with
distinct centers, ${\bf x}$ and ${\bar {\bf x}_{{\bf u}- {\{\bf
x}\}}}$, which implies that
 ${\cal
V}_{{\bf u} - \{{\bf x}\}}$ is perpendicular to the direction
vector of the line
$$
{L}_{{\bf u},{\bf x}}(\alpha)={\bar {\bf x}_{{\bf u} - \{{\bf
x}\}}} + \alpha \left({\bf x}-{\bar {\bf x}_{{\bf u} - \{{\bf
x}\}}} \right) \quad (\alpha \in \openR)
$$
passing through these centers and that the points of ${\bf u} -
\{{\bf x}\}$ are equidistant from each point of $L_{{\bf u},{\bf
x}}$. In particular, all points of ${\bf u}- \{{\bf x}\}$ are
equidistant from $L_{{\bf u},{\bf x}}(1/\n)={\bar {\bf x}_{\bf
u}}$,   hence so are all points of ${\bf u}$. Because ${\bf
x}-{\bar {\bf x}_{{\bf u} - \{{\bf x}\}}} \perp {\cal V}_{{\bf u}
- \{{\bf x}\}}$ but ${\bf x}-{\bar {\bf x}_{{\bf u} - \{{\bf
x}\}}} \in {\cal V}_{{\bf u}}$, $\mbox{dim}\left( {\cal V}_{\bf u}
\right)=\mbox{dim}\left( {\cal V}_{{\bf u}- \{{\bf x}\}}
\right)+1$. By orthogonality $||{\bar {\bf x}_{\bf u}} - {\bar
{\bf x}_{{\bf u}- \{{\bf
x}\}}}||^2=r_{\sigma,\n}^2-r_{\sigma,\n-1}^2$ and does not depend
on ${\bf x}$. Using the fact that ${\bar {\bf x}_{{\bf u}- \{{\bf
x}\}}},{\bar {\bf x}_{\bf u}},$ and ${\bf x}$ all lie on $L_{{\bf
u},{\bf x}}$ in tandem with orthogonality gives
$2\sigma^2=r_{\sigma,\n-1}^2+(s_{\sigma,\n}+r_{\sigma,\n})^2$;
solving these two equations for $r_{\sigma,\n}$ and
$s_{\sigma,\n}$ finishes the induction.
 \end{proof}

 \vspace{3mm}

\noindent {\it Proof of the theorem.}\ \  Let ${\bf X}= ({\bf
x}_1,\ldots,{\bf x}_n)$, an element of $\openR^{p \times n}$.
Since ${\bf S}_{{\bf T}({\bf u})}={\bf S}_{\bf u}$ for any
translation ${\bf T}$, we can assume without loss of generality
that the members of ${\bf u}$ have already been centered by
subtraction of their mean, so ${\bar {\bf x}}_{{\bf u}}={\bf 0}$
and in general letting ${\bf B}_{\bf v}:=|{\bf v}|{\bf S}_{\bf v}$
we have
\bea
\label{null} {\bf B}_{\bf u}=\sum_{{\bf x} \in {\bf u}} {\bf x}
{\bf x}'= {\bf X}{\bf X}'.
\ena
Assuming that the points are equidistant, we infer from
(\ref{null}) and the lemma that
\begin{eqnarray*}
{\bf B}_{\bf u}{\bf x}&=&\sum_{{\bf y} \in {\bf u}- \{{\bf
x}\}}{\bf y}{\bf y}'{\bf x}+{\bf x}{\bf x}'{\bf x} =
(r_{\sigma,n}^2-\sigma^2)\sum_{{\bf y} \in {\bf u}- \{{\bf
x}\}}{\bf y} + r_{\sigma,n}^2{\bf x}\\
&=&   (\sigma^2-r_{\sigma,n}^2){\bf
x}+ r_{\sigma,n}^2{\bf x} = \sigma^2 {\bf x}
\end{eqnarray*}
for each ${\bf x}$ in ${\bf u}$.  Hence ${\bf B}_{\bf u}{\bf
x}=\sigma^2 {\bf I}_p{\bf x}$ on ${\cal
V}_{\bf u}$.  Since $\mbox{dim}({\cal V}_{\bf u})=p$ by the lemma,
${\bf B}_{\bf u}=\sigma^2 {\bf I}_p$.

For the converse, assume that ${\bf B}_{\bf u}=\sigma^2{\bf I}_p$.
Note that the matrix
$$
{\bf A}=\sigma^{-2}{\bf X}'{\bf X} \in {\bf R}^{n \times n}
$$
is symmetric, ${\bf A}'={\bf A}$, and idempotent, ${\bf
A}^2=\sigma^{-4}{\bf X}'{\bf X} {\bf X}'{\bf X}=\sigma^{-4}{\bf
X}'{\bf B}_{\bf u}{\bf X}={\bf A}$. Hence ${\bf A}$ is an
orthogonal projection, and therefore has rank equal to its trace,
$$
\mbox{rank}({\bf A})=\mbox{tr}({\bf A})=\sigma^{-2}\mbox{tr}({\bf
X}'{\bf X}) = \sigma^{-2}\mbox{tr}({\bf X}{\bf X}')=
\sigma^{-2}\mbox{tr}({\bf B}_{\bf u})=p,
$$
using the cyclic invariance of the trace. With ${\bf 1}_n \in {\bf
R}^n$ the vector with all components equal to 1, ${\bf A}{\bf
1}_n={\bf 0}$ by virtue of ${\bar {\bf x}}_{{\bf u}}=0$. By the
rank plus nullity theorem the null space of ${\bf A}$ has
dimension one, and must therefore equal $\mbox{span}({\bf 1}_n)$,
the span of ${\bf 1}_n$. Hence ${\bf A}={\bf I}_n-\frac{1}{n}{\bf
1}_n{\bf 1}_n'$, as this is the unique orthogonal projection of
rank $p$ with null space $\mbox{span}({\bf 1}_n)$. As the entries
of ${\bf A}$ are $\sigma^{-2}$ times the inner products of the
vectors in ${\bf u}$, the squared interpoint distances between
${\bf x}_i \not = {\bf x}_j$ equals
$$
||{\bf x}_i-{\bf x}_j||^2 =2\left({\bf x}_i'{\bf x}_i-{\bf
x}_i{\bf x}_j \right) =2\sigma^2 \left(
\left(1-\frac{1}{n}\right)+\left(
\frac{1}{n}\right)\right)=2\sigma^2.
$$
$\qed$

We remark that once the matrix ${\bf A}$ is determined to have
constant off-diagonal entries, the proof may also be completed by
induction in the following more geometric way: Assume that $n
> 2$, the base case being trivial. Any $p$ points in $\openR^p$
lie in a hyperplane of dimension $p-1$, and for ${\bf x} \in {\bf
u}$ let ${\cal H}$ denote the hyperplane which contains ${\bf u} -
\{{\bf x}\}$, the space ${\cal V}_{{\bf u} - \{{\bf x}\}}$
translated by ${\bar {\bf x}}_{{\bf u} - \{{\bf x}\}}$. The inner
products ${\bf x}'{\bf y}$ for all ${\bf y} \in {\bf u}-\{{\bf
x}\}$, being the off-diagonal elements of ${\bf A}$, are equal,
and therefore, ${\bf x}'{\bf y}={\bf x}'{\bar {\bf x}}_{{\bf u} -
\{{\bf x}\}}$, so ${\bf x} \perp {\bf y}-{\bar {\bf x}}_{{\bf u} -
\{{\bf x}\}}$. Hence ${\bf x} \perp {\cal V}_{{\bf u} - \{{\bf
x}\}}$, and since ${\bar {\bf x}}_{{\bf u} - \{{\bf x}\}}=-{\bf
x}/p$, we conclude ${\bf x} - {\bar {\bf x}}_{{\bf u} - \{{\bf
x}\}} \perp {\cal V}_{{\bf u} - \{{\bf x}\}}={\cal H}-{\bar {\bf
x}}_{{\bf u} - \{{\bf x}\}}$.

Now let ${\bf T}$ be the translation ${\bf T}{\bf y}={\bf y}-{\bar
{\bf x}_{{\bf u} - \{{\bf x}\}}}$, and, with $\{{\bf e}_i\}_{1 \le
i \le p}$ the standard basis, ${\bf O}$ the rotation that maps
${\bf T}{\bf x}$ to $\beta e_p$ where $\beta=||{\bf x}-{\bar {\bf
x}_{{\bf u} - \{{\bf x}\}}}||$. That is, ${\bf V}({\bf
x})=\beta{\bf e}_p$ for ${\bf V}={\bf O}{\bf T}$, and
\bea
\label{BVu} &&{\bf B}_{{\bf V}({\bf u}- \{{\bf x}\})}+\beta^2{\bf
e}_p{\bf e}_p'={\bf O}\left({\bf B}_{{\bf T}({\bf u}- \{{\bf x})\}}+
{\bf T}{\bf x}({\bf T}{\bf x})'\right){\bf O}'\\
\nonumber &=& {\bf O}\left({\bf B}_{{\bf u}- \{{\bf x}\}}+{\bf
T}{\bf x}({\bf T}{\bf x})'\right){\bf O}'={\bf O}{\bf B}_{{\bf
T}({\bf u})}{\bf O}' ={\bf O}{\bf B}_{{\bf u}}{\bf
O}'=\sigma^2{\bf O} {\bf I}_p {\bf O}'= \sigma^2 {\bf I}_p.
\ena
Since ${\cal V}_{{\bf u} - \{{\bf x}\}} \perp {\bf x}-{\bar {\bf
x}_{{\bf u} - \{{\bf x}\}}}$, ${\bf V}({\cal H}) \subset
\openR^{p-1} \times \{0\}$, and we can consider the points ${\bf
V}({\bf u}- \{{\bf x}\})$ as lying in $\openR^{p-1}$. By
(\ref{BVu}), the $(p-1) \times (p-1)$ submatrix $[{\bf B}_{{\bf
V}({\bf u} - \{{\bf x}\})}]_{1 \le i,j \le p-1}$ equals
$\sigma^2{\bf I}_{p-1}$, so applying the induction hypotheses to
${\bf V}({\bf u}- \{{\bf x}\})$ we conclude that the interpoint
distances of ${\bf u}- \{{\bf x}\}$, unchanged by ${\bf V}$, are
all $2\sigma^2$. The induction is completed by noting that this is
true for each ${\bf x}$ in ${\bf u}$.

{\em Acknowledgement} The authors thank Richard E. Stone for
bringing a shortcoming in the original version of this work to
their attention.

\begin{tabular}{ll}
Ian Abramson & Larry Goldstein\\
Department of Mathematics  & Department of Mathematics\\
University of California, San Diego& University of Southern California\\
La Jolla, CA 92093-0112 & Los Angeles, CA  90089-2532\\
iabramson@ucsd.edu & larry@math.usc.edu
\end{tabular}


\begin{thebibliography}{1}

\bibitem[1]{1}
G. A. F. Seber and A. J. Lee  {\em Linear Regression Analysis},
John Wiley, New York, 2003.
\end{thebibliography}
\end{document}